\documentclass[12pt,a4paper]{article}

\usepackage{a4wide,amssymb,amsmath,amsthm,xspace,epsfig}

\begin{document}

\newcommand{\N}{\mathbb{N}}
\newcommand{\R}{\mathbb{R}}
\newcommand{\Z}{\mathbb{Z}}
\newcommand{\Q}{\mathbb{Q}}
\newcommand{\C}{\mathbb{C}}
\newcommand{\PP}{\mathbb{P}}

\newcommand{\LL}{\Bbb L}

\newcommand{\esp}{\vskip .3cm \noindent}
\mathchardef\flat="115B

\newcommand{\lev}{\text{\rm Lev}}

\def\ut#1{$\underline{\text{#1}}$}
\def\CC#1{${\cal C}^{#1}$}
\def\h#1{\hat #1}
\def\t#1{\tilde #1}
\def\wt#1{\widetilde{#1}}
\def\wh#1{\widehat{#1}}
\def\wb#1{\overline{#1}}

\def\restrict#1{\bigr|_{#1}}

\newtheorem{lemma}{Lemma}[section]

\newtheorem{thm}[lemma]{Theorem}

\newtheorem{defi}[lemma]{Definition}
\newtheorem{conj}[lemma]{Conjecture}
\newtheorem{cor}[lemma]{Corollary}
\newtheorem{prop}[lemma]{Proposition}
\newtheorem*{prob}{Problem}
\newtheorem{ques}[lemma]{Question}
\newtheorem*{rem}{Remark}
\newtheorem{exemples}[lemma]{Examples}
\newtheorem{exemple}[lemma]{Example}

\title{Notes on the od-Lindel\"of property}
\date{\today}
\author{Mathieu Baillif}
\maketitle

\abstract{A space is od-compact (resp. od-Lindel\"of) provided any cover by open dense sets has a finite (resp. countable) subcover.
We first show with simple examples that these properties behave quite poorly under finite or countable unions. 
We then investigate the relations between Lindel\"ofness, od-Lindel\"ofness and linear Lindel\"ofness (and similar relations with `compact'). 
We prove in particular that if a $T_1$ space is od-compact, then the subset of its non-isolated points is compact.
If a $T_1$ space is od-Lindel\"of, we only get that the subset of its non-isolated points is linearly Lindel\"of. Though, 
Lindel\"ofness follows if the
space is moreover locally openly Lindel\"of (i.e. each point has an open Lindel\"of neighborhood).}

\section*{Note}
After the completion (and publication) of this paper, the author became aware that 
stronger results were obtained
by Mills and E. Wattel some time ago in \cite{MillsWattel}. See also \cite{Blair:1983}.
Some of our results were actually proved even before in \cite{Katetov:1947}.

\section{Introduction}

In the middle of an argument involving Baire theorem, we noticed that we did not need the space under scrutinity to be really Lindel\"of, but rather that any cover of it 
by open {\em dense} sets
had a countable subcover. We then wondered whether this alternative definition of Lindel\"ofness, called here {\em od-Lindel\"ofness}, 
was interesting in itself, as well as the similarly defined notion of od-compactness.
These notes
are the results of our musings, which may be summarized as follows.

$\bullet$ od-compact spaces behave quite horribly when taking unions, even when just two subspaces are involved, and there are even completely
          metrizable spaces
          that behave bad in this respect. A finite union of od-compact {\em closed} spaces is od-compact, though. 
          On the other hand a countable union of od-Lindel\"of closed spaces does not need to be od-Lindel\"of.

\smallskip
$\bullet$ The image of an od-compact space under a continuous map is not always od-compact, and the same holds for od-Lindel\"of spaces.
          However the properties are preserved when the map is open. Moreover, the image of a $T_1$ od-compact space by a closed map is od-compact.

\smallskip
$\bullet$ Trivial examples of od-compact spaces are the discrete ones. But in a way they are the only non-compact ones.
          In fact, the subset of non-isolated points of a $T_1$ od-compact space is compact.
          For od-Lindel\"ofness, our results are not that strong.
          First, an od-Lindel\"of $T_1$ space that does not contain a clopen uncountable discrete subset {\em and} which is locally openly Lindel\"of 
          is Lindel\"of (see below for undefined terminology).
          If one drops the last assumption, then we could only obtain that
          the space is linearly Lindel\"of. 
          (In fact, the result for od-compact spaces follows from the equivalence of the linarly compact and compact notions.)
          Moreover, the examples we know of linearly Lindel\"of spaces that are {\em not} Lindel\"of happen to be non-od-Lindel\"of as well.

\medskip
We have not found older references to these od-notions, but since our examples and proofs are rather elementary, we would not be
surprised if some of our results already appeared somewhere. 
Perhaps the above points provide an explanation for this absence in the literature: 
od-compact and od-Lindel\"of properties are not `robust' at 
all, and moreover (at least for the compact case), differ only slightly from the usual compact and Lindel\"of notions. 
However, we would be interested in finding a non-trivial example of od-Lindel\"of non-Lindel\"of space, or in showing that there
is none.

\smallskip
This note is organized as follows.
In Section \ref{sec:defs} we give the definitions and show some equivalences. 
In Section \ref{sec:unions} we investigate the behavior of the od- properties when taking unions.
Then, we prove the above mentioned theorem relating od-Lindel\"ofness with Lindel\"ofness  in Section \ref{sec:od} while
the relation with linear Lindel\"ofness is shown in Section \ref{sec:odlin}.
The short Section \ref{sec:images} contains the above mentioned results about images of od-Lindel\"of and od-compact spaces under open and closed maps.
We included a short appendix containing classical results featuring compactness, Lindel\"ofness and complete accumulation points.

Most of this note does not contain or use technicalities beyond the basics of topology and elementary ordinal/cardinal manipulation, and is
fairly self contained.
However, some of the examples we shall give are classical spaces of set-theoretic topology for which we will just
give a reference,
and we shall have a few words about more recent constructions of linearly Lindel\"of non-Lindel\"of spaces. We shall refer to 
the articles where these spaces were described for more details.


\section{Definitions}\label{sec:defs}

`Space' always means `topological space'. 
We use the greek letters $\alpha,\beta,\gamma$ for ordinals, and $\kappa,\lambda,\tau$ for cardinals.
We denote by $\wb{B}$ and $int(B)$ the closure and interior of a subset of a space.

\begin{defi} Let $X$ be a space.
   \begin{itemize} 
   \item[$\bullet$]
     $L(X)$ is the smallest cardinal $\kappa$ such any open cover of $X$ has a subcover of cardinality $<\kappa$.
     $X$ is compact if $L(X)\le\omega$ and Lindel\"of if $L(X)\le\omega_1$, and more generally  Lindel\"of$_\kappa$ if $L(X)\le\kappa$.
   \item[$\bullet$]
     $\ell L(X)$ is the smallest cardinal $\kappa$ such any open cover of $X$, which is a chain for the inclusion 
     (in short: a chain-cover), has a subcover of cardinality $<\kappa$.
     $X$ is linearly compact if $\ell L(X)\le\omega$ and linearly Lindel\"of if $\ell L(X)\le\omega_1$, 
     and more generally linearly Lindel\"of$_\kappa$ if $\ell L(X)\le\kappa$.
   \item[$\bullet$]
     $odL(X)$ is the smallest cardinal $\kappa$ such any cover of $X$ by open dense sets (in short: an od-cover) has a subcover of cardinality $<\kappa$.
     $X$ is od-compact if $odL(X)\le\omega$ and od-Lindel\"of if $odL(X)\le\omega_1$,
     and more generally  od-Lindel\"of$_\kappa$ if $odL(X)\le\kappa$.
   \end{itemize}
\end{defi}

Beware than in a lot of texts, the similar Lindel\"of degree of a space is defined a bit differently (for instance, $L(\R)=\omega_1$, while its 
Lindel\"of degree is $\omega$).
We chosed this definition because it seems to enable shorter statements when compact spaces are also involved.
Of course, Lindel\"of$_\omega$ and Lindel\"of$_{\omega_1}$ are synonyms of compact and Lindel\"of.
Notice also that we do {\em not} assume any separation axiom for compactness and Lindel\"ofness, though it is not difficult to show
that one can assume our spaces to be $T_0$ by taking Kolmogorov quotients.
It was shown long ago that linearly compact spaces are compact, see the appendix.

\begin{exemples}
\quad \\
$\bullet$ Any Lindel\"of$_\kappa$ space is od-Lindel\"of$_\kappa$ and linearly Lindel\"of$_\kappa$.\\
$\bullet$ Any space with the discrete topology is od-compact (in fact, $odL(X)=2$). 
\end{exemples}

Recall the following elementary lemma:

\begin{lemma}\label{Lclosed}\ \\
   a) For a topological space $X$ and $\kappa\ge\omega$, $L(X)=\kappa$ iff $L(Y)=\kappa$ for each closed $Y\subset X$,
   iff given a family of closed sets with empty intersection, there is a subfamily of cardinality $<\kappa$ with empty intersection.
   \\
   b) If $X$ is a union of $\kappa$ spaces $X_\alpha$ with $L(X_\alpha)\le\lambda$ for a regular $\lambda$, then $L(X)\le \kappa\cdot\lambda$.
\end{lemma}

When od- properties are concerned, we obtain:

\begin{lemma}\label{ODclosed}
   Are equivalent:\\
   a) $odL(X) \le \kappa$,\\   
   b) Any cover of $X$ by open sets such that at least one is dense has a subcover of cardinality $<\kappa$.\\
   c) $odL(Y) \le \kappa$ for each closed $Y\subset X$,\\
   d) $L(Y)\le\kappa$ for each $Y\subset X$ closed and nowhere dense.
\end{lemma}

In particular: a space is od-compact (resp. od-Lindel\"of) iff each of its closed nowhere dense subsets is compact (resp. Lindel\"of).

\proof \ \\
   a) and b) are easily seen equivalent: given an open cover $U_\alpha$ for $\alpha\in\lambda$ with $U_0$ dense, then
   the sets $U_0\cup U_\alpha$ form an od-cover.\\
   a) $\Rightarrow$ c) If $C$ is closed in $X$ and $U$ is an open dense set in $C$, then there is a $V$ open in $X$ with $C\cap V=U$, and
   $V\cup(X- C)$ is dense.\\
   c) $\Rightarrow$ a) Immediate.\\
   a) $\Rightarrow$ d) If $L(Y) > \kappa$ for some nowhere dense closed $Y$, then 
     given a cover of $Y$ witnessing this fact we find an od-cover of $X$ taking the union of each member with $X-Y$.\\
   d) $\Rightarrow$ b) Given an open cover $U_\alpha$ ($\alpha\in\lambda$) of $X$ such that $U_0$ is dense, set $B_\alpha=X-U_\alpha$. Then 
      $B_0$ is nowhere dense, and $\cap_{\alpha\in\lambda}B_\alpha = \cap_{\alpha\in\lambda}(B_0\cap B_\alpha) = \varnothing$. 
      Since $L(B_0)\le\kappa$, there is a subfamily of the $B_\alpha$ of cardinality $<\kappa$ with empty intersection by Lemma \ref{Lclosed} a).
      The corresponding family of $U_\alpha$ cover $X$.
\endproof


\section{Unions of od-Lindel\"of$_\kappa$ spaces}\label{sec:unions}

The od-covering properties behave in a quite horrible manner when taking unions.

\begin{exemple}\label{horribleunion1}
  For each cardinal $\kappa\ge\omega$, there is a $T_1$ space $X$ with $odL(X)=\kappa^+$, which satisfies
  $X=X_0\sqcup X_1$, where $X_0$ is compact and $X_1$ closed and discrete (so $odL(X_0)=\omega$, $odL(X_1)=2$).
\end{exemple}
If $\kappa=\omega$ set $\gamma =\omega\cdot\omega$, otherwise set $\gamma=\kappa$.
$X$ is given by $(\gamma+1)\times\{0,1\}$ with the following topology.
Let the topology on $X_0=(\gamma+1)\times\{0\}$ be the usual order topology of $\gamma+1$, $X_0$ is thus compact.
The neighborhoods of $(\alpha,1)$ are all the subsets of $X$ than can be written as 
$(U-F)\times\{0\}\sqcup F\times\{1\}$, where $U$ is open in $\gamma+1$, and $F\subset U$ is a finite set containing $\alpha$.
Then $X_1=(\gamma+1)\times\{1\}$ is discrete in $X$.
One shows easily that $X$ is $T_1$ (but not Hausdorff).
Set $U$ to be the open set given by $X_0$ union $\{(\alpha,1)\,:\,\alpha\text{ successor}\}$. (Recall that $\{\alpha\}$ is open in $\gamma+1$
iff $\alpha$ is successor.) $U$ is dense in $X$.
For each limit $\alpha\in\gamma+1$, set $U_{\alpha}=U\cup({\alpha}\times\{1\})$, $U_\alpha$ is then open and dense. 
The od-cover by the $U_\alpha$s does not have any subcover of cardinality $<\kappa$.

\smallskip
The same type of idea can be used to obtain:

\begin{exemple}\label{horribleunion2}
  For each cardinal $\kappa\ge\omega$, there is a completely metrizable space $X$ with $odL(X)\ge\kappa^+$, which satisfies
  $X=X_0\sqcup X_1$, where $X_0,X_1$ are discrete and $X_0$ is closed (so $odL(X_0)=odL(X_1)=2$).
\end{exemple}
Take $X$ to be a disjoint union of clopen copies $J_\alpha$ ($\alpha\in\kappa$) of $\{0\}\cup\{\frac{1}{m}\,:\,m\in\omega\}$, each with its usual topology.
A complete metric on $X$ is given by the usual distance for two points in the same $J_\alpha$, while two points in two different $J_\alpha$ are at distance $2$.
Then $X_0$ is the union of the $0$ points, while $X_1$ is its complement.
The od-cover given by the $U_\alpha$ defined as $X_1\cup J_\alpha$ has no proper subcover.

\smallskip
Still another example in the same vein, this time for (non)-od-Lindel\"of spaces,
showing that we cannot even trust a subspace of a (non-metrizable) $2$-manifold:

\begin{exemple}\label{horribleunion3}
  There is a subspace $S$ of a $2$-manifold with $odL(S)=L(S)=(2^\omega)^+$, such that $S=A\sqcup B$, with 
  $A$ closed discrete, and $B\simeq\R^2$ (so $odL(A)=2$, $odL(B)=\omega_1$).
\end{exemple}
This example is the 
subset of the (separable version of the)
Pr\"ufer surface, with $A$ being given by taking one point in each boundary component, and $B$ is the interior (i.e. 
the surface minus the boundary components). See for instance the appendix in \cite{Spivak:vol1} for a description.
The idea is essentially a `manifold equivalent' to the tangent disk topology on the half plane which is described
in \cite[Example 82]{CEIT}. Both contain a closed nowhere dense discrete subset of cardinality $2^\omega$, and are thus 
non-od-Lindel\"of$_{2^{\omega}}$ by Lemma \ref{ODclosed}.

\medskip
Examples \ref{horribleunion1} to \ref{horribleunion3} all make use of a closed discrete subset whose complement is dense.
It is easy to see that one cannot hope to find two closed sets whose union behaves that bad:

\begin{lemma}\label{unionclosed} 
  Let $\kappa$ be an infinite cardinal.
  If $X=X_0\cup\cdots\cup X_n$ is a finite union of closed od-Lindel\"of$_\kappa$ subsets for $i=1,\dots,n$, then $odL(X)\le\kappa$.
\end{lemma}
\proof
  We prove it for two subsets, the general case follows by induction. Let thus $X=X_0\cup X_1$,
  and let $B\subset X$ be closed and nowhere dense. 
  We shall show that $L(B\cap X_i)\le\kappa$ for $i=0,1$, which implies $L(B)\le\kappa$ and the result
  by Lemma \ref{ODclosed}. We may thus assume first that $B\subset X_0$, the other case being entirely symmetric. 

\begin{figure}[h]
    \centering
    \epsfig{figure=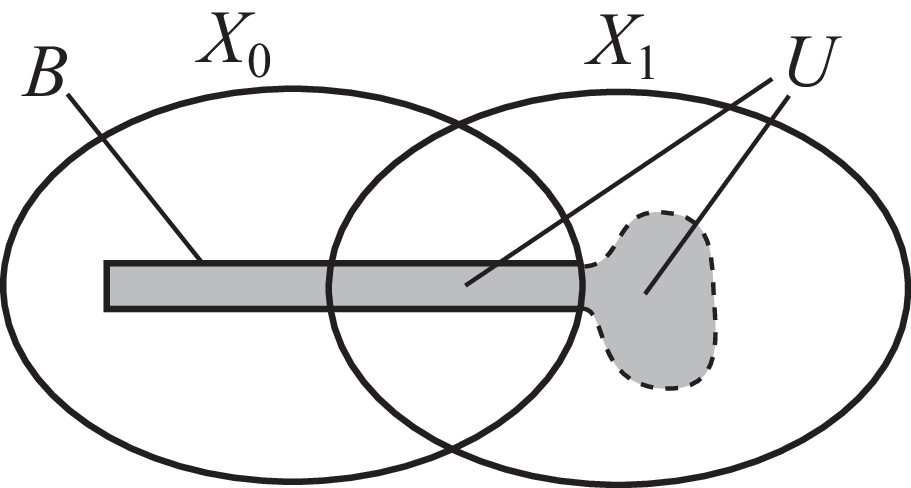, width=50mm}
    \caption{\label{fig:1} Proof of Lemma \ref{unionclosed}.}
\end{figure}

  Denote by $int_0$ the interior for the induced topology in 
  $X_0$. 
  If $int_0(B)$ is empty, then $B$ is nowhere dense in $X_0$ and $L(B)\le\kappa$ by Lemma \ref{ODclosed}.
  If not, let $U\subset X$ be open with
  $U\cap X_0=int_0(B)$ (as in Figure \ref{fig:1}).
  Notice that $L(B - int_0(B))\le\kappa$, since $B-int_0(B)$ is closed and nowhere dense in $X_0$.
  If $(U\cap X_0) - X_1\not=\varnothing$, then 
  $(U-X_1)\subset U\cap X_0 \subset B$, so $B$ is not nowhere dense in $X$. 
  Thus $(U\cap X_0)\subset X_1$, so 
  $$
    U=(U\cap X_0)\cup  (U\cap X_1) \subset X_1,
  $$ 
  $U\cap X_1$ is open in $X$, and contains $int_0(B)$.
  It follows that $int_0(B)$ is nowhere dense in $X_1$
  (otherwise for some $W$ open in $U$ and thus in $X$, $W\subset int_0(B)\subset B$),
  so $L(\wb{int_0(B)})\le\kappa$, where the closure is taken in $X_1$ (or in $X$ since $X_1$ is closed).
  Since $B = \wb{int_0(B)}\cup \bigl(B - int_0(B)\bigr)$,
  $L(B)\le\kappa$.
\endproof

Also, there is no Hausdorff space that has the properties of Example \ref{horribleunion1}:

\begin{lemma}\label{unioncompact} \ \\
  a) If $X=X_0\cup X_1$ with $X_0$ closed, $L(X_0)\le\kappa$, and $odL(X_1)\le\kappa$, then $odL(X)\le\kappa$.\\
  b) If $X=X_0\cup X_1$ is Hausdorff, with $X_0$ compact and $odL(X_1)\le\kappa$, then $odL(X)\le\kappa$.
\end{lemma}
\proof
  a) 
  Take $B\subset X$ to be nowhere dense and closed. 
  By Lemma \ref{Lclosed}, $L(B\cap X_0)\le\kappa$, and since $X_1-X_0$ is open, 
  $B\cap (X_1-X_0)$ is nowhere dense in $X_1$, so $L(B\cap (X_1-X_0))\le\kappa$.
  Thus $L(B)\le\kappa$. \\
  b) Since $X$ is Hausdorff, $X_0$ is closed, and we apply a). 
\endproof

The situation with countable unions is bad even for $\sigma$-discrete (i.e. a countable union of closed discrete subspaces) spaces.

\begin{exemple}
  There are Tychonov locally compact $\sigma$-discrete non-od-Lindel\"of spaces.
\end{exemple}
Such a space is thus a countable union of closed od-compact subspaces but is non-od-Lindel\"of.
Any Hausdorff Aronszajn special $\omega_1$-tree $T$ 
with the order topology is such an example, since it is a countable union of antichains which are closed discrete subspaces and thus 
od-compact
(see for instance \cite{Nyikos:trees} for definitions, especially Theorem 4.11).
Moreover such a tree is locally compact and Hausdorff, and thus Tychonov. 
However, if one denotes the members of $T$ at height $\alpha$ by $T_\alpha$ and the set of limit ordinals by $\Lambda$,
the od-cover given by $U_\alpha = T - {\displaystyle \cup_{\beta\in\Lambda,\, \beta>\alpha}}T_\alpha$ has no countable subcover.


\section{od-Lindel\"ofness in locally (openly) Lindel\"of spaces}\label{sec:od}

From now on, `cardinal' means `infinite cardinal'.
There are various definitions of local Lindel\"ofness in the literature. We opted for the following terminology for clarity.

\begin{defi} Let $\tau$ be a regular cardinal. 
          A space $X$ is locally [openly] Lindel\"of$_\tau$ provided each of its points possesses a Lindel\"of$_\tau$ neighborhood [which is open].
\end{defi}

Recall that the notions agree for regular spaces (and regular cardinals $\tau\ge\omega_1$):
\begin{lemma}\label{OLandL}
   Let $\tau\ge\omega_1$ be regular, and $X$ be a regular space.
   Then
   $X$ is locally Lindel\"of$_\tau$ if and only if it is locally openly Lindel\"of$_\tau$, if and only if it has a basis of closed Lindel\"of neighborhoods.
\end{lemma}

A proof can be found for instance instance by combining \cite[Theorem 2.3]{DissanayakeSastry} and
\cite[Prop. 1.1]{Balogh:2002} (the result is stated for $\tau=\omega_1$, but the proof works in general).
When the space is not regular, the result does not hold anymore:

\begin{exemples}\ \\
   $\bullet$
   The everywhere doubled line (see \cite{BaillifGabard}) is a locally Euclidean $T_1$ space 
   which is locally (openly) Lindel\"of but does not have a basis of closed Lindel\"of neighborhoods. \\
   $\bullet$ The half disk topology (Example 78 in \cite{CEIT}) is a Hausdorff example of such a space.
\end{exemples}

(Neither example is od-Lindel\"of, though.)
The goal of this section is the proof of the following theorem.

\begin{thm}\label{thm1}
   Let $\tau$ be a regular cardinal, $X$ be a $T_1$  
   locally openly Lindel\"of$_\tau$ space with $odL(X)\le\tau$. Then either $L(X)\le\tau$, or 
   there is a clopen discrete subset of cardinality $\ge\tau$ in $X$ . 
\end{thm}

(Note that when $\kappa=\omega$, Theorem \ref{thmoLlL0} c) below is much stronger.)
The core of the proof is essentially contained in the next lemma.

\begin{lemma}\label{cruciallemma}
   Let $\tau,\lambda$ be regular cardinals, and let $X$ be a $T_1$ space with $odL(X)\le\tau$.
   Let $Y\subset X$ be closed, and $Z_\alpha$ be open for $\alpha\in\lambda$, such that
   $Y\subset \cup_{\alpha\in\lambda}Z_\alpha$, 
   $Z_\alpha\subset Z_\beta$ whenever $\alpha<\beta<\lambda$, $\wb{Z_\alpha}\not\supset Y$
   for each $\alpha$. 
   Then, either $\lambda<\tau$, or 
   $\lambda\ge\tau$ and there is a discrete subset $D\subset Y$, clopen in $X$, of cardinality $\ge\lambda$.
\end{lemma}

\proof 
   We shall define $x_\alpha\in Y$ and $f:\lambda\to\lambda$ as follows. Set $f(0)=0$. 
   Given $f(\alpha)$, choose $x_\alpha\in Y - \wb{Z_{f(\alpha)}}$, and set $f(\alpha+1)$ to be the smallest $\beta$ such that
   $Z_\beta\ni x_\alpha$. When $\alpha$ is limit, set $f(\alpha)=\sup_{\beta<\alpha}f(\alpha)$.
   Since $\lambda$ is regular, $f(\alpha)$ and $x_\alpha$ are defined for each $\alpha<\lambda$.
   Set $U_\alpha=Z_{f(\alpha)}$, the $U_\alpha$s have the same properties as the $Z_\alpha$s, and 
   \begin{equation}\label{eq:one}
      x_\alpha\in U_{\alpha+1}-\wb{U_\alpha}.
   \end{equation} 
   Let $E$ be the set of $\alpha$ such that $(\wb{U_\alpha}-U_\alpha)\cap Y\not=\varnothing$.
   If $E$ is cofinal in $\lambda$, then letting 
   $$
     V_\alpha = X - \cup_{\alpha<\beta<\lambda}(\wb{U_\beta}-U_\beta)\cap Y,
   $$
   we get a cover of $X$ by open dense subsets without any subcover of cardinality $<\lambda$, which implies $\lambda<\tau$.
   The $V_\alpha$s are indeed open, since any point $y$ in the closure of $\cup_{\alpha<\beta<\lambda}(\wb{U_\beta}-U_\beta)\cap Y$
   must be in $Y$, and thus in $U_\gamma$ for some $\gamma$ which we can take minimal, $\gamma$ is therefore successor and equal to some
   $\xi+1$, with $\xi>\alpha$. Then $y\in\wb{U_\xi}$, otherwise $U_\gamma - \wb{U_\xi}$ is a neighborhood of $y$ that intersects no
   $(\wb{U_\beta}-U_\beta)$, 
   and $y\not\in U_\xi$ by minimality of $\gamma$. So $y\in (\wb{U_\xi}-U_\xi)\cap Y$.
   
   We may thus assume that $E\subset\alpha<\lambda$ for some $\alpha$, 
   and in fact that $E=\varnothing$.
   Set $B=\{x_\alpha\,:\alpha\in\lambda\}$. 
   By (\ref{eq:one}), $\{x_\alpha\}$ is open in $\wb{B}$ for the induced topology.
   Since $X$ is $T_1$ and the $U_\alpha$ cover $Y$ which is closed, $\wb{B}-B$ is contained in the union
   of the $(\wb{U_\alpha}-U_\alpha)\cap Y$ for limit $\alpha$, which are empty, thus $B$ is closed, as well as any of its subsets.
   If the interior of ${B}$ is contained in some $U_\alpha$, then ${B}-U_\alpha$ is closed and nowhere dense, so by Lemma \ref{ODclosed},
   $L({B}-U_\alpha)\le\tau$. But the $U_\beta$ for $\beta<\lambda$ cover it and there is no subcover of cardinality $<\lambda$ by
   (\ref{eq:one}), and thus $\lambda<\tau$.
   So let us assume now that $int(B)$ is not contained in any $U_\alpha$. 
   Then the $\alpha$ for which $\{x_\alpha\}$ is open in $X$ are cofinal in $\lambda$. 
   Call $D$ the union of all these open $\{x_\alpha\}$, then $D$ is clopen and $|D|=\lambda$.
\endproof

Another auxiliary result that we shall use:

\begin{lemma}\label{lemma1thm1}
  If there is a subset $U\subset X$ which is open, Lindel\"of$_\tau$ and such that $\wb{U}$ is not Lindel\"of$_\tau$,
  then $odL(X)>\tau$.
\end{lemma}
\proof
  Otherwise, the nowhere dense closed subset $\wb{U}-U$ would be Lindel\"of$_\tau$ by Lemma \ref{ODclosed}, and $\wb{U}=(\wb{U}-U)\cup U$ as well.
\endproof

We now start the proof of Theorem \ref{thm1} in earnest.

\proof[Proof of Theorem \ref{thm1}]
  Suppose that $L(X)>\tau$.
  By Lemma \ref{lemma1thm1}, we can assume that 
  \begin{equation}\label{eq1}
      \wb{U} \text{ is Lindel\"of$_\tau$ whenever $U$ is open and Lindel\"of$_\tau$.}
  \end{equation}
  We will build open subsets $X_\alpha$ for ordinals $\alpha$.
  Let $X_0\subset X$ be any open Lindel\"of$_\tau$ subset, and build $X_\alpha$
  as follows. If $\alpha$ is limit, take $X_\alpha=\cup_{\beta<\alpha}X_\beta$.
  If $\alpha=\beta+1$ and $\wb{X_\beta}-X_\beta\not=\varnothing$, 
  take a Lindel\"of$_\tau$ open neighborood $U_x$ of each 
  $x\in\wb{X_\beta}-X_\beta$.
  If $\wb{X_\beta}$ is Lindel\"of$_\tau$ (and thus, $\wb{X_\beta}-X_\beta$ as well), 
  extract a subcover 
  $U_{x_i}$ ($i\in\tau_0<\tau$), and set $X_\alpha=X_\beta\cup(\cup_{i\in\tau_0}U_{x_i})$
  ($X_\alpha$ is then Lindel\"of$_\tau$ by Lemma \ref{Lclosed}). 
  If such a subcover does not exist, set
  $X_\alpha=X_\beta\cup(\cup_{x\in\wb{X_\beta}-X_\beta}U_x)$.
  If $\wb{X_\beta}=X_\beta \not= X$, choose an open Lindel\"of$_\tau$ set $U$ disjoint from
  $X_\beta$, and set $X_\alpha= X_\beta\cup U$.

  By construction, we have $\wb{X_\beta}\subset X_\alpha$ whenever $\beta < \alpha$.
  For some $\alpha$, $X=X_\alpha$. Take $\alpha$ to be minimal
  with this property.
  Let $\beta$ be the supremum of 
  $\{\gamma < \alpha\,:\,\wb{X_\gamma}\text{ is Lindel\"of$_\tau$}\}$.
  Then $\wb{X_{\beta}}$ is not Lindel\"of$_\tau$, otherwise by construction and (\ref{eq1}), so would be $\wb{X_{\beta+1}}$. Likewise,
  $X_\beta$ is not Lindel\"of$_\tau$.
  If $\beta$ is successor, $\wb{X_{\beta -1}}$ would be Lindel\"of$_\tau$, so $X_\beta$ as well, and 
  similarly, if $\beta$ is limit with $cf(\beta)<\tau$, $X_\beta=\cup_{\gamma<\beta}\wb{X_\gamma}$ would be a union 
  of less than $\tau$ Lindel\"of$_\tau$ spaces, and therefore Lindel\"of$_\tau$ by Lemma \ref{Lclosed}. 
  Thus, $cf(\beta)\ge\tau$. 
  We now have two cases. (The case $\beta = \alpha$ is contained in the first one, with $V=\varnothing$.)
  \begin{itemize}
  \item[1)] There is an open $V\supset (\wb{X_\beta}-X_\beta)$ such that
     the set $\{\gamma<\beta\,:\, (X_\beta - (V\cup X_\gamma)) \not=\varnothing\}$ is 
     cofinal in $\beta$.
     Then, $X$ satisfies the assumptions of Lemma \ref{cruciallemma} with $Y=\wb{X_\beta} - V$ and
     $\lambda=cf(\beta)\ge\tau$, and $X$ contains a clopen
     discrete subset of cardinality $\ge\tau$.
     
  \item[2)] For any open set $V\supset (\wb{X_\beta}-X_\beta)$, there is a 
     $\gamma<\beta$ such that 
     $(X_{\beta}-X_\gamma)\subset V$.
  
  Suppose that $\wb{X_\beta} - X_\beta$ is Lindel\"of$_\tau$, and 
  let $\langle U_i\,:\,i\in I\rangle$ be an open cover of $\wb{X_\beta}$.
  Extract a subcover of $\wb{X_\beta} - X_\beta$ of cardinality $<\tau$, and choose $\gamma < \beta$ such
  that $X_\gamma$ is Lindel\"of$_\tau$ and $(X_{\beta}-X_\gamma)$ is included in the union of this
  subcover. Adding a  subcover of $\wb{X_\gamma}$ of cardinality $<\tau$ and putting everything together
  yields a subcover of $\wb{X_\beta}$ of the same cardinality, so 
  $\wb{X_\beta}$ is Lindel\"of$_\tau$, and $\wb{X_{\beta+1}}$ as well, contradicting the definition
  of $\beta$. Therefore $\wb{X_\beta} - X_\beta$ is not Lindel\"of$_\tau$.

  Let thus $U_i$ ($i\in I$) be a cover of $\wb{X_\beta} - X_\beta$ without
  subcover of cardinality $<\tau$. Set $W_i= U_i \cup X_\beta \cup (X-\wb{X_\beta})$, which yields a cover
  of $X$ by open dense sets with the same property, a contradiction since $X$ is od-Lindel\"of$_\tau$.
  \end{itemize}
\endproof

In view of the impressive list given in \cite{Gauld:Met}, it might be interesting to notice the following corollary:

\begin{cor}
  A manifold is metrizable if and only if it is od-Lindel\"of.
\end{cor}

\proof
  A manifold is metrizable iff all its connected components are metrizable, so we may assume the manifold to be connected.
  A manifold is locally compact and its singletons are not open, so
  it cannot possess an open discrete subset, hence is Lindel\"of if and only if od-Lindel\"of. 
  We conclude by recalling that Lindel\"ofness and metrizability are equivalent for connected manifolds.
\endproof

The next lemma yields more consequences of Theorem \ref{thm1}:

\begin{lemma}\label{isolatedlemma}
   Let $\tau$ be a regular cardinal, $X$ an od-Lindel\"of$_\tau$ space, $D$ the subspace of its isolated points. 
   Then $X-D$ does not contain a clopen discrete subset of cardinality $\ge\tau$.
\end{lemma}
\proof
  Notice that $D$ is open and discrete, so
  $X-D$, being closed, is od-Lindel\"of$_\tau$ by Lemma \ref{ODclosed}.
  Suppose that $X-D$ contains a clopen (in $X-D$) discrete subset $D_0$ of cardinality $\ge\tau$.
  Then $D\cup\{x\}$ is a neighborhood of $x$ for each $x\in D_0$, and setting
  $V_x=\{x\}\cup(X-D_0)$ yields an od-cover without subcover of cardinality $<\tau$.
\endproof

It follows immediately:

\begin{cor}\label{thm1c}
  Let $\tau$ be a regular cardinal,
  $X$ be a locally openly Lindel\"of$_\tau$ space with
  $odL(X)\le\tau$. Let $D\subset X$ be the subset of isolated points.
  Then $X-D$ is Lindel\"of$_\tau$.
\end{cor}

We shall later  relax the local openly Lindel\"ofness assumption, so let us introduce a notation.

\begin{defi}
  Let $\tau>\omega$ be a regular cardinal, and
  $X$ be a topological space.
  $$
    \begin{array}{rcl}
       \mathsf{L}_\tau(X)&=&\{x\in X\,:\,\exists\text{ an open Lindel\"of$_\tau$ }U\ni x\}\\
        \mathsf{NL}_\tau(X)&=&X- \mathsf{L}_\tau(X).
    \end{array}
  $$
  We denote by $\mathsf{C}(X)$ the subset containing the points possessing a compact neighborhood, and set $\mathsf{NC}(X)=X-\mathsf{C}(X)$.
\end{defi}

It is immediate from the definition that $\mathsf{L}_\tau(X)$ and $\mathsf{C}(X)$ are open. 
There are simple spaces with $\mathsf{NC}(X)$ (resp. $\mathsf{NL}_\tau(X)$) consisting of just one point:
the cone $[0,1]\times Y / \ (0,y)\sim(0,z)$ over any locally compact (resp. locally openly Lindel\"of$_\tau$) $Y$ which is not
compact (resp. Lindel\"of$_\tau$). 

\begin{thm}\label{thm1b} 
  Let $\kappa\ge\omega_1$ be a regular cardinal and
      $X$ be a $T_1$ space such that $odL(X)\le\kappa$ and $L(\mathsf{NL}_\kappa(X))\le\kappa$. Then, either $L(X)\le\kappa$, or $X$
      contains a clopen discrete subset of cardinality $\ge\kappa$.
\end{thm}

\proof
Notice that if $L(X)\le\kappa$, then $\mathsf{NL}_\kappa(X)=\varnothing$.  
We have two cases. 
\begin{itemize}
   \item[i)] There is some open $U\supset \mathsf{NL}_\kappa(X)$ such that $L(X-U)>\kappa$.
    \\
   We repeat the proof of Theorem \ref{thm1} in $X-U$ (which is od-Lindel\"of$_\kappa$)
   and apply Lemma \ref{cruciallemma} in case 1) for $Y=(\wb{X_\beta} - V)\cap (X-U)$, yielding the same result.
  \item[ii)] For all open $U\supset \mathsf{NL}_\kappa(X)$, $L(X-U)\le\kappa$.
   \\
   In this case, $L(X)$ will be $\le\kappa$. Indeed, given a cover of $X$ by $V_i$, $i\in I$,
   let $V_{i_k}$ for $k\in J$ be a subcover of $\mathsf{NL}_\kappa(X)$ of cardinality $<\kappa$. 
   Then, $X-\cup_{k\in J}V_{i_k}$ being Lindel\"of$_\kappa$,
   is covered by $<\kappa$ many more $V_i$.
\end{itemize}
\endproof


\section{od- and linear-Lindel\"ofness}\label{sec:odlin}

Here, we show the relations between od- and linear-Lindel\"ofness.
First, an easy theorem.

\begin{thm}\label{thmoLlL0} 
   Let $\kappa$ be a regular cardinal. \\
   a) The subspace of non-isolated points of a $T_1$ od-Lindel\"of$_\kappa$ space is linearly Lindel\"of$_\kappa$.\\
   b) If the subspace of non-isolated points of a space is Lindel\"of$_\kappa$, the space is od-Lindel\"of$_\kappa$.\\
   c) A $T_1$ space is od-compact iff the subspace of its non-isolated points is compact.
\end{thm}
\proof
  a) Let $D$ contain the isolated points of $X$ and set $Z=X-D$. Then by Lemma \ref{ODclosed} $odL(Z)\le\kappa$, and
     $Z$ does not have a clopen discrete subset of cardinality $\ge\kappa$ by Lemma \ref{isolatedlemma}.
     Let $U_\alpha$ ($\alpha\in\lambda$) be 
     a chain-cover of $Z$.
     We may assume $\lambda$ to be regular.
     If some $U_\alpha$ is dense in $Z$, then 
     each $U_\beta$ for $\beta\ge\alpha$ is such, so there is a subcover of 
     $Z$ of cardinality $<\kappa$. 
     We may now assume that none of the $U_\alpha$ is dense in $Z$.
     But then $X$ satisfies the hypotheses of Lemma \ref{cruciallemma} for $Y=Z$,
     which yields $\lambda<\kappa$.\\
  b) By Lemma \ref{unioncompact} a).\\
  c) By Corollary \ref{lincomcom} below, a linearly compact space is compact, the result follows thus from a) and b).
\endproof

When $\kappa>\omega$, one can get a finer result (though not as good as in the compact case):

\begin{thm}\label{thmoLlL} 
      Let $X$ be a $T_1$ space with $odL(X)\le\kappa$ for a regular $\kappa\ge\omega_1$, and $D\subset X$ be the subset of isolated points.
      Then $X=D\sqcup X_0\sqcup X_1$, where $X_0\cup D$ is open, $L(\wb{X_0})\le\kappa$, $X_1$ is closed, $\ell L(X_1)\le\kappa$, and
      any open set $U$ with $U\cap X_1\not=\varnothing$ satisfies $L(U)>\kappa$.\\
\end{thm}
 
\proof
     Set $Z=X-D$, again $odL(Z)\le\kappa$, and
     $Z$ does not have a clopen discrete subset of cardinality $\ge\kappa$.
     Set $X_0=\mathsf{L}_\kappa(Z)$, $X_1=\mathsf{NL}_\kappa(Z)$.
     By Lemma \ref{ODclosed}, $odL(\wb{X_0})\le\kappa$ and $odL(X_1)\le\kappa$. 
     Notice that $\mathsf{NL}_\kappa(\wb{X_0})=\wb{X_0}\cap X_1$ is closed and nowhere dense, so $L(\mathsf{NL}_\kappa(\wb{X_0}))\le\kappa$,
     and by Theorem \ref{thm1b}, $L(\wb{X_0})\le\kappa$. 

     We now repeat the proof of Theorem \ref{thmoLlL0}. Let $U_\alpha$ ($\alpha\in\lambda$) be 
     a chain-cover of $X_1$.
     As above
     we may assume that none of the $U_\alpha$ is dense in $X_1$.
     But then $Z$ satisfies the hypotheses of Lemma \ref{cruciallemma} for $Y=\wb{X_1}\supset\cup_{\alpha\in\lambda}U_\alpha$,
     which yields again $\lambda<\kappa$.
\endproof

\begin{exemples} There are linearly Lindel\"of non-od-Lindel\"of spaces.
\end{exemples}
These spaces are examples of linearly Lindel\"of non-Lindel\"of (abbreviated $\ell$LnL below) spaces found in the literature, which happen 
to be non-od-Lindel\"of.
\begin{itemize}    
  \item[$\bullet$] Probably the first example of an $\ell$LnL space was given by Mi\v s\v cenko in \cite{Mishenko}. It is a Tychonoff space, defined as 
    the subset of $R=\Pi_{i\in\omega}(\omega_i + 1)$ by the 
    union $\cup_{k\in\omega}R_k$ with $R_k= \bigl(\Pi_{i=0,\dots,k-1}(\omega_i + 1)\bigr)\times \bigl(\Pi_{i\in\omega,\,i\ge k}\omega_i\bigr)$.
    (As usual, we denote by $\omega_i$ the $i$-th cardinal above $\omega = \omega_0$, and by $\omega_\omega$ the sup of these $\omega_i$.)
    The proof given in \cite{Mishenko} can be easily adapted 
    to show that the od-cover given by the $\Gamma_{\alpha,i}$, defined for $i\in\omega$ and $\alpha\in\omega_\omega$ as
    the subset of points whose $i$-th coordinate is not a limit ordinal $\ge\alpha$, does not admit a subcover of cardinality $<\aleph_\omega$, so 
    this space is not od-Lindel\"of. 
  \item[$\bullet$]
    Arhangel'skii and Buzyakova \cite[Example 4.1]{ArhanBuzy:1999} gave a description of another Tychonoff $\ell$LnL space $X$,
    which is a subspace of 
    $\mathcal{D}^A$, where $\mathcal{D}$ is the discrete space $\{0,1\}$ and $A$ is discrete with cardinality $\aleph_\omega$.
    $X$ is the subspace consisting of the points that have less than $\aleph_\omega$ coordinates equal to $1$.
    They show that $X$ is pseudocompact since it contains a dense countably compact subspace, and non-compact since 
    it is not closed in $\mathcal{D}^A$. 
    It happens that $X$ is non-od-Lindel\"of. Indeed, fix an uncountable $A_0\subset A$ such that $|A-A_0|=\aleph_\omega$, 
    and let $B$ be the subset of $X$ consisting of points whose coordinates in $A_0$ are all $0$.
    Then $B$ is closed and nowhere dense (since it does not contain a basic open set, where only a finite number of coordinates are fixed).
    But $B$ is homeomorphic to $X$, and thus non Lindel\"of, so by Lemma \ref{ODclosed} $X$ is non-od-Lindel\"of.
    The modified version in \cite{ArhanBuzy:1998} has the same property.
  \item[$\bullet$]
    Kunen \cite{kunen:2002, kunen:2005} found locally compact $\ell$LnL spaces.
    Recall that a locally compact space is Tychonoff and thus regular, so by Lemma \ref{OLandL}, $X$ is locally openly Lindel\"of, thus
    $\mathsf{NL}_{\omega_1}(X)=\varnothing$. 
    A linearly Lindel\"of space does not contain an uncountable clopen discrete subset, so by Theorem
    \ref{thm1} $X$ is not od-Lindel\"of.
\end{itemize}

These results and examples raise the following questions:

\begin{ques}
   Is there a $T_1$ space which does not contain a clopen uncountable discrete subset that is od-Lindel\"of and non-Lindel\"of~?
\end{ques}

\begin{ques}
   What conditions should be added to linear Lindel\"ofness to ensure that a space is od-Lindel\"of~?
\end{ques}


\section{Images of od-Lindel\"of spaces}\label{sec:images}

Notice that the continuous image of an od-compact space may be violently non-od-compact:

\begin{exemple} 
  Denote by $\kappa_d$ the cardinal $\kappa$ with the discrete topology, while $\kappa$ is endowed with the usual order topology.
  Then $odL(\kappa_d)= 2$, while $odL(\kappa)=cf(\kappa)$, and the identity map
  $\kappa_d\to\kappa$ is continuous.
\end{exemple}

However we have preservation if the map is open, and also if the map is closed and $X$ is $T_1$ and od-compact.
The proof of the latter fact uses Theorem \ref{thmoLlL0} c). We found neither an easier proof (which we believe should exist) 
nor a general result for od-Lindel\"of$_\kappa$ spaces with uncountable $\kappa$.

\begin{lemma}\label{openclosed}
   Let $X$, $Y$ be spaces, and $f:X\to Y$ be continuous. \\
   a) If $f$ is open then $odL(f(X))=odL(X)$.\\
   b) If $f$ is closed and $X$ is $T_1$ and od-compact, then $f(X)$ is od-compact.
\end{lemma}

\proof 
     In both cases we may assume that $f(X)=Y$.

      a) First, $f(X)$ is open in $Y$, so a relatively open subset of $f(X)$ is indeed open.
      Let $\{U_j\,:\,j\in J\}$ be an od-cover of $f(X)$.
      If $f^{-1}(U_j)$ misses some open nonempty $W\subset X$, then $U_j\cap f(W)=\varnothing$, which is impossible. Thus
      $\{f^{-1}(U_j)\,:\,j\in J\}$ is an od-cover, and we conclude by extracting a subcover and mapping it through $f$.

      b) Let $D$ be the set of isolated points of $X$, then $X-D$ is closed and compact by Theorem \ref{thmoLlL0} c) 
      so $f(X-D)$ is closed and compact as well.
      We now show that the points in $Y-f(X-D)$ are isolated, by Lemma \ref{unioncompact} a) this yields that $f(X)$ is od-compact.
      Let $x\in D$ be such that $f(x)\notin f(X-D)$. Define the open subset 
      $$
        Z_x=\{z\in D\,:\, f(z)=f(x)\},
      $$
      then $\{f(x)\} = Y-f(X-Z_x)$ is open, which shows that $f(x)$ is isolated.
\endproof


\section{Appendix: Classical results on linearly Lindel\"of and compact spaces}\label{sec:classical}

Here we recall some classical basic results, due to Alexandroff and Urysohn \cite{AlexandroffUrysohn}.
Consider the following properties for a space $X$ and a regular infinite cardinal $\kappa$:
\begin{align}\label{align1}\tag{CAP$(\kappa)$}
   &  \begin{array}{l}
        \text{If $B$ is a subset of regular cardinality $\ge\kappa$, it has a point of complete}\\
        \text{accumulation.} 
      \end{array}
    \\
      \label{align2}\tag{CAP$^+(\kappa)$}
    & \begin{array}{l}
        \text{If $B$ is a subset of cardinality $\ge\kappa$, it has a point of complete}\\
        \text{accumulation.} 
      \end{array}
\end{align}

Then, we have:

\begin{lemma}\label{capcap+} $X$ satisfies 
   CAP$(\omega)$ iff  it satisfies CAP$^+(\omega)$.
\end{lemma}

\proof  
  CAP$^+(\omega)$ implies trivially CAP$(\omega)$, we thus show the other implication.
  Let $\kappa\ge\omega$ be minimal such that 
  there is some $B\subset X$ with $|B|=\kappa$ without complete accumulation point, $\kappa$ must be singular and $>\omega$ by CAP$(\omega)$.
  Thus, for all infinite $\lambda < \kappa'$, there is an accumulation point $x_\lambda$ of $B$ such that any open set
  containing $x_\lambda$ intersects $B$ in at least $\lambda$ points.
  Let $\tau=cf(\kappa)<\kappa$, and
  Let $f:\tau\to\kappa$ be a cofinal map.
  Since $C=\{x_{f(\alpha)}\,:\,\alpha\in\tau\}$ has a cardinality less than $\kappa$ but $\ge\omega$, it possesses a complete 
  accumulation point $x$. 
  (In this part of the proof we really need $\omega$.)
  Thus, any open $U\ni x$ contains $x_\lambda$ for a subset of $\lambda$ cofinal in $\kappa$.
  Hence, it intersects $B$ in more than $\lambda$ points for each $\lambda<\kappa$, and therefore in $\kappa$ points.
\endproof

\begin{thm}\label{cpccap} Let $X$ be a space and $\kappa$ be regular.\\
   a) $X$ satisfies CAP$^+(\omega)$ iff $L(X)=\omega$ (i.e. $X$ is compact).\\
   b) If $X$ satisfies CAP$^+(\kappa)$ then $L(X)\le\kappa$.\\
   c) $X$ satisfies CAP$(\kappa)$ iff $\ell L(X)\le\kappa$.\\
\end{thm}

\proof
  a)  
  Assume $X$ to be compact, and let $B\subset X$ be infinite. If there is no complete accumulation point for $B$, then
  for each $x\in X$ there is an open set $U_x\ni x$ with $|U_x\cap B|<|B|$. Taking a finite subcover, this yields that
  $|B|$ is a finite sum of smaller cardinals, which is impossible. The converse is included in b).

  b) Let $\kappa$ be regular.
  Suppose that $L(X)>\kappa$, and let $\kappa'$ be minimal such that there exists an open cover 
  $\langle U_\alpha\,:\,\alpha\in\kappa'\rangle$ of $X$ without a subcover of cardinality $<\kappa$. Set $V_\alpha = \cup_{\beta<\alpha}U_\alpha$.
  If for some $\alpha<\kappa'$ we have $V_\alpha = X - E$, with $|E|<\kappa'$, then letting $\beta(x)$ be the smallest
  $\beta$ such that $x\in U_\beta$, we get that 
  $$
  \langle U_\beta\,:\, \beta<\alpha \text{ or }\beta=\beta(x)\text{ for some } x\in E\rangle
  $$
  is a cover of $X$ by less than $\kappa'$ open sets, thus by minimality of $\kappa'$ there is a cover of cardinality $<\kappa$, a contradiction.
  Thus, for each $\alpha$ there is $x_\alpha\not\in V_\alpha$. Hence $x_\alpha\not\in U_\beta$ for each $\beta<\alpha$, 
  and $B=\{x_\alpha\,:\,\alpha\in\kappa'\}$ has no complete accumulation point. (Because each $x\in X$ belongs to some $U_\beta$ which contains
  $<\kappa'$ points of $B$.)

  c)  
  Assume that $\ell L(X)\le\kappa$, and let $B=\{x_\alpha\,:\,\alpha <\kappa'\}$ for some regular $\kappa'\ge\kappa$.
  Set $B_\beta = \{ x_\alpha\,:\,\beta\le\alpha <\kappa'\}$, and 
  $U_\beta = X-\wb{B_\beta}$.
  Then $\langle U_\beta\,:\,\beta\in\kappa'\rangle$ is a chain for the inclusion.
  If it covers $X$, we may extract a subcover of cardinal $<\kappa$, and since $\kappa'$ is regular, there is some $\beta<\kappa'$ (the sup of the
  indices in the subcover) with $U_\beta = X$. Thus $B_\alpha$ is empty for each $\alpha>\beta$, a contradiction.
  Therefore there is some $x\in X$ such that $x\not\in U_\beta$ (that is,
  $x\in\wb{B_\beta}$) for all $\beta$. Given an open set $U\ni x$, for each $\beta$ there is an $\alpha\ge\beta$ with
  $x_\alpha\in U$. The regularity of $\kappa'$ implies then that $|U\cap B|=\kappa'$, so $x$ is a complete accumulation point.

  Conversely, given an open cover $\langle U_j\,:\,j\in J\rangle$ of $X$ which is a chain
  and does not possess a subcover of cardinality $<\kappa$, let $\lambda$ be minimal such that
  there is a cofinal map $f:\lambda\to J$. Then $\lambda$ is regular, and writing $V_\alpha$ for $U_{f(\alpha)}$, 
  $\langle V_\alpha\,:\,\alpha\in \lambda\rangle$ is a cover of $X$, which does not possess a subcover of cardinality $<\kappa$.
  For each $\alpha\in\lambda$ let $x_\alpha\not\in V_\alpha$, then 
  $B=\{x_\alpha\,:\,\alpha\in\lambda\}$ has no complete accumulation point, 
  because each $V_\alpha$ contains less than $\lambda$ points of $B$, and they cover $X$. This contradicts CAP$(\kappa)$.
\endproof

Notice that the last part of the proof does not work if one takes a cover that is not a chain.
Moreover, the converse implication of b) does not hold: $\omega_\omega$ induced with the order topology is Lindel\"of but 
it does not posses a point of complete accumulation.

The corollary we used in Theorem \ref{thmoLlL0} follows immediately:
\begin{cor}\label{lincomcom}
  A space is compact iff it is linearly compact.
\end{cor}


\begin{thebibliography}{10}

\bibitem{AlexandroffUrysohn}
P.~Alexandroff and P.~Urysohn.
\newblock {\em {M\'emoire} sur les espaces topologiques compacts}.
\newblock Number~14 in Proceedings of the section of mathematical sciences.
  Koninklijke Nederlandse Akademie van Wetenschappen te Amsterdam, 1929.

\bibitem{ArhanBuzy:1999}
A.V. ArhangelÕskii and R.~Z. Buzyakova.
\newblock On linearly {Lindel\"of} and strongly discretely {Lindel\"of} spaces.
\newblock {\em Proc. Amer. Math. Soc.}, 127(8):2449--2458, 1999.

\bibitem{ArhanBuzy:1998}
A.V. ArhangelÕskii and R.Z. Buzyakova.
\newblock Convergence in compacta and linear {Lindel\"ofness}.
\newblock {\em Comment. Math. Univ. Carolin.}, 39(1):159--166, 1998.

\bibitem{BaillifGabard}
M.~Baillif and A.~Gabard.
\newblock Manifolds: Hausdorffness versus homogeneity.
\newblock {\em Proc. Amer. Math. Soc.}, 136:1105--1111, 2008.

\bibitem{Balogh:2002}
Z.~Balogh.
\newblock Locally nice spaces and axiom {R}.
\newblock {\em Topology Appl.}, 125(2):335--341, 2002.

\bibitem{Blair:1983}
R.L. Blair.
\newblock Some nowhere densely generated topological properties.
\newblock {\em Comm. Math. Univ. Carolinae}, 24(3):465--479, 1983.

\bibitem{DissanayakeSastry}
U.~N.~B. Dissanayake and K.~P.~R. Sastry.
\newblock Locally {Lindel\"of} spaces.
\newblock {\em Indian J. pure appl. Math.}, 18(10):876--881, 1987.

\bibitem{Katetov:1947}
M.~Kat\v etov.
\newblock On the equivalence of certain types of extensions of topological
  spaces.
\newblock {\em \v Casepis p\v est. mat. fys.}, 72:101--106, 1947.

\bibitem{Gauld:Met}
D.~Gauld.
\newblock Metrisability of manifolds.
\newblock Preprint arXiv:0910.0885v1.

\bibitem{kunen:2002}
K.~Kunen.
\newblock Locally compact linearly {Lindel\"of} spaces.
\newblock {\em Comment. Math. Univ. Carolin.}, 43(1):155--158, 2002.

\bibitem{kunen:2005}
K.~Kunen.
\newblock Small locally compact linearly {Lindel\"of} spaces.
\newblock {\em Topology Proc.}, 29(1):193--198, 2005.

\bibitem{MillsWattel}
C.F. Mills and E.~Wattel.
\newblock {\em Nowhere densely generated topological properties}, pages
  191--198.
\newblock Number 116 in Math. Centre Tracts. Math. Centrum, Amsterdam, 1979.

\bibitem{Mishenko}
A.~S. {Mi\v s\v cenko}.
\newblock Finally compact spaces.
\newblock {\em Soviet Math. Dokl.}, 145:1199--1202, 1962.

\bibitem{Nyikos:trees}
P.~Nyikos.
\newblock Various topologies on trees.
\newblock In P.R. Misra and M.~Rajagopalan, editors, {\em Proceedings of the
  Tennessee Topology Conference}, pages 167--198. World Scientific, 1997.

\bibitem{Spivak:vol1}
M.~Spivak.
\newblock {\em A comprehensive introduction to differential geometry},
  volume~I.
\newblock Publish or Perish, Wilmington, 1979.
\newblock Second edition.

\bibitem{CEIT}
L.~A. Steen and J.~A.~Seebach Jr.
\newblock {\em Counter examples in topology}.
\newblock Springer Verlag, New York, 1978.

\end{thebibliography}
\end{document}